\documentstyle{amsppt}
\NoRunningHeads \magnification=\magstephalf
\hcorrection{.5truein}
\NoBlackBoxes \refstyle{C} \loadbold

\def\fett{
\textfont0=\tenbf \scriptfont0=\sevenbf \scriptscriptfont0=\fivebf
\textfont1=\tenbi \scriptfont1=\sevenbi \scriptscriptfont1=\fivebi
\textfont2=\tenbsy \scriptfont2=\sevenbsy
\scriptscriptfont2=\fivebsy \rm
\bf}

\font\teneuf=eufm10 
\font\seveneuf=eufm7 
\font\fiveeuf=eufm5
\newfam\euffam \def\euf{\fam\euffam\teneuf}
\textfont\euffam=\teneuf \scriptfont\euffam=\seveneuf
\scriptscriptfont\euffam=\fiveeuf %
\def\Afr{{\euf A}}

\let\fett=\hbox
\let\cal=\Cal

    \overfullrule=0pt

\def\thfill{\null\nobreak\hfill}
\def\endbeweis{\thfill\vbox{\hrule
     \hbox{\vrule\hbox to 5pt{\vbox to 5pt{\vfil}\hfil}\vrule}\hrule}}

\def\ref#1{\par\noindent\hangindent3\parindent
\hbox to 3\parindent{#1\hfil}\ignorespaces}

\topmatter

\title Existence de filtrations admissibles sur des isocristaux
\endtitle

     \author
J.-M. Fontaine \footnote{{Institut Universitaire de France et UMR
8628
du CNRS, Math\'{e}matique, Universit\'e de Paris-Sud, B\^atiment
425,
91405 Orsay Cedex, France}\hfill\ }
     et M. Rapoport
\footnote{ Mathematisches Institut der Universit\"at zu K\"oln,
Weyertal
86-90, D--50931 K\"oln, Deutschland}
     \endauthor

\endtopmatter

\define\Dbf{\bold D}
\define\G{\bold G}
\define\N{\bold N}
\define\Q{\bold Q}
\define\R{\bold R}
\define\Z{\bold Z}
\document

\bigskip
\noindent {\bf 0 -- Enonc\'{e} des r\'{e}sultats}
\medskip\noindent

Soit $F$ une extension finie de $\Q_p$. Soit $k$ un corps parfait de
caract\'{e}ristique $p$ contenant le corps r\'{e}siduel $k_F$ de $F$.
Notons $W(k)$ (resp. $W(k_F)$) l'anneau des vecteurs de Witt
\`{a}
coefficients dans $k$ (resp. $k_F$) et posons
$L=F\otimes_{W(k_F)}W(k)$.
C'est un corps, extension finie totalement ramifi\'{e}e du corps des
fractions de $W(k)$. Soient $q$ le cardinal de $k_F$, $\sigma_0$
l'automorphisme de $W(k)$ induit par fonctorialit\'{e} par
l'automorphisme $x\mapsto x^{q}$ sur $k$ et $\sigma$
l'automorphisme
$\roman{id}\otimes\sigma_0$ de $L$. Un
{\it isocristal} ({\it relativement \`{a}} $(F,k)$) est un
espace vectoriel
$D$ de dimension finie sur $L$,
muni d'un endomorphisme $\varphi$ bijectif et $\sigma$-lin\'eaire.
Soit
$(\Q^d)_+=\ \{ (\nu_1,\nu_2,\ldots,\nu_d)\in \Q^d\mid\nu_1\geq
\nu_2\geq
\ldots \geq \nu_d\}$. A un isocristal $(D,\varphi)$ de dimension
$d$
on associe son vecteur de Newton (= la suite des pentes en ordre
{\it
d\'ecroissant}), $\nu(D,\varphi)\in (\Q^d)_+$. Soit
$(\Z^d)_+=(\Q^d)_+\cap
\Z^d$. Soit $\mu=(\mu_1,\ldots, \mu_d)\in (\Z^d)_+$. Soit  ${\Cal
F}^\bullet$ une
filtration de $D$ par des sous-$L$-espaces vectoriels,
d\'ecroissante, exhaustive et s\'epar\'{e}e, index\'{e}e par $\Z$~; on
dit que ${\Cal F}^{\bullet}$
{\it est  de type $\mu$} si $${\roman{dim}}\  gr_i^{{\Cal F}^\bullet}
(D)=
\# \{ j;\ \mu_j=i\}\ \ .\leqno(0.1)$$ Autrement dit, les sauts de la
filtration ${\Cal F}^\bullet$ sont les $\mu_j$ ($j=1,\ldots, d)$ et la
taille du saut en $\mu_j$ est donn\'ee par la multiplicit\'e de $\mu_j$
dans $\mu$. Toute filtration ${\Cal F}^\bullet$ a un type $\mu({\Cal
F}^\bullet)\in (\Z^d)_+$ bien d\'etermin\'e.

Une filtration
${\Cal F}^\bullet$ de l'isocristal $(D,\varphi)$ de dimension $d$
est
dite {\it admissible} si les deux conditions suivantes sont
satisfaites~:

(i) on a $$\vert\mu({\Cal F}^\bullet)\vert = \vert\nu (D,\varphi)\vert \ \
$$ (on a not\'e ici $\vert\nu\vert =\sum_{i=1}^d\nu_i$ pour
$\nu=(\nu_1,\ldots, \nu_d)\in\Q^d$)~;

(ii) soit $(D', \varphi')$ un sous-isocristal de $(D,\varphi)$ et soit
${\Cal F}^{\prime\bullet}$ la filtration induite par ${\Cal F}^\bullet$ sur
$D'$~; on a
$$\vert\mu({\Cal F}^{\prime\bullet})\vert\leq \vert\nu(D', \varphi')\vert\
\ .$$

\medskip\noindent {\bf Remarques~:} 1) Pour tout
$\lambda=(\lambda_1,\lambda_2,\ldots,\lambda_d)\in(\Q^{d})_+$,
notons
$P(\lambda)$ le {\it polygone associ\'e \`{a} $\lambda$}, i.e. le
polygone convexe (autrement dit, \`{a} pentes croissantes) du plan
r\'{e}el d'origine $(0,0)$ dont les pentes
sont les $\lambda_j$, la longueur de la projection du segment de
pente
$\lambda_j$ sur l'axe des $x$ \'{e}tant \'{e}gale \`{a} la
multiplicit\'{e} de $\lambda_j$ dans $\lambda$.
Le nombre rationnel $\vert\nu(D,\varphi)\vert$,
not\'{e} $t_N(D)$ dans [F], est
toujours un entier~; les points de coordonn\'{e}es $(0,0)$ et
$(d,\vert\nu(D,\varphi)\vert)$
sont les extr\'{e}mit\'{e}s du polygone $P(\nu(D,\varphi))$ (parfois
appel\'{e} {\it polygone de Newton} de l'isocristal
$(D,\varphi)$). De m\^{e}me, $\vert\mu({\Cal F}^\bullet)\vert$ est
not\'{e} $t_H({\Cal F}^{\bullet})$ dans [F]~; les points de
coordonn\'{e}es $(0,0)$ et
$(d,\vert\mu({\Cal F}^\bullet)\vert)$
sont les extr\'{e}mit\'{e}s du polygone $P(\mu({\Cal F}^{\bullet}))$
(parfois appel\'{e} {\it polygone de Hodge}
de la filtration).

2) Dans le cas o\`{u} $F=\Q_p$, dire qu'une filtration est admissible
signifie qu'elle est {\it faiblement admissible} au sens de [F]. Le
r\'{e}sultat principal de [CF] signifie qu'elle est alors \'{e}galement
admissible au sens de [F].  Si $\bar L$ d\'{e}signe une cl\^{o}ture
alg\'{e}brique de $L$, on dispose donc ({\it loc.cit.)} d'une
\'{e}quivalence naturelle entre la cat\'{e}gorie des isocristaux
(relativement \`{a} $(\Q_p,k)$) munis d'une filtration admissible, et la
cat\'{e}gorie des repr\'{e}sentations $p$-adiques cristal\-lines de
${\roman{Gal}}(\bar L/L)$.

3) Ce r\'{e}sultat s'\'{e}tend au cas o\`{u} $F$ est une extension finie
arbitraire de $\Q_p$. Soit $C_p$ le compl\'{e}t\'{e} de $\bar L$ pour la
topologie $p$-adique. Appelons {\it $F$-repr\'{e}sentation de}
$\roman{Gal}(\bar L/L)$ la donn\'{e}e d'un $F$-espace vectoriel de
dimension finie $V$ muni d'une action lin\'{e}aire et continue de
$\roman{Gal}(\bar L/L)$. Pour une telle repr\'{e}sentation, notons
$V_{C_p,F}$ le noyau de l'application $C_p$-lin\'{e}aire naturelle de
$C_p\otimes_{\Q_p}V$ dans $C_p\otimes_F V$. On dispose alors d'une
\'{e}quivalence naturelle entre la cat\'{e}gorie des isocristaux
(relativement \`{a} $(F,k)$) munis d'une filtration admissible et la
sous-cat\'{e}gorie pleine de la cat\'{e}gorie des $F$-repr\'{e}sentations
de $\roman{Gal}(\bar L/L)$ dont les objets sont les $V$ qui sont
cristallines en tant que repr\'{e}sentations $p$-adiques et v\'{e}rifient
la condition $$(\ast)\hbox{ le }C_p\hbox{-espace vectoriel
}V_{C_p,F}\hbox{ est engendr\'{e} par les \'{e}l\'{e}ments fixes par }
\roman{Gal}(\bar L/L)\ \ .$$ Lorsque $V$ provient d'un groupe formel
$\Phi$, cette derni\`{e}re condition signifie que $\Phi$ est un
$O_F$-module formel (cf. par exemple, [De]). Par exemple, soit $\pi$ une
uniformisante de $F$ et soit $(D,\varphi)=(L,\pi\sigma)$~; la seule
filtration admissible est celle pour laquelle $gr_1^{{\Cal
F}^\bullet}(D)\not=0$. La $F$-repr\'{e}sentation de $\roman{Gal}(\bar
L/L)$ associ\'{e}e est de dimension $1$ et c'est la duale de la
restriction \`{a} $\roman{Gal}(\bar L/L)$ de celle que fournit un groupe
formel de Lubin-Tate pour $F$ correspondant \`{a} $\pi$. Soit
$\eta:\roman{Gal}(\bar L/L)\to F^{\times}$ le caract\`{e}re qui d\'{e}finit
l'action du groupe de Galois. Si $\tau$ est un $\Q_p$-automorphisme non
trivial de $F$, la $F$-repr\'{e}sentation de dimension $1$ d\'{e}finie par
$\tau\eta$ est encore cristalline mais elle ne v\'{e}rifie pas $(\ast)$.

    Les remarques (2) et (3) ne seront pas utilis\'{e}es dans la suite.
\medskip\noindent

Sur $(\Q^d)_+$, on dispose de l'ordre partiel  pour lequel
$\lambda\leq\lambda'$
si $\lambda_1+\lambda_2\ldots+\lambda_r\leq
\lambda'_1+\lambda'_2\ldots+\lambda'_r$, pour $r=1,2,\ldots,d-1$ et
$\lambda_1+\lambda_2\ldots+\lambda_d=
\lambda'_1+\lambda'_2\ldots+\lambda'_d$. Dire que
$\lambda\leq\lambda'$ \'{e}quivaut \`{a} dire que $P(\lambda)$ est
au-dessus de
$P(\lambda')$ et que ces deux polygones ont m\^{e}mes
extr\'{e}mit\'{e}s.

\medskip\noindent
{\bf Th\'eor\`eme 1.}  {\it Supposons $k$ alg\'{e}briquement clos.
Soit $(D,\varphi)$ un isocristal de
dimension $d$. Soit $\mu\in(\Z^d)_+$. Pour qu'il existe une filtration
faiblement admissible de type $\mu$ sur $D$ il faut et il suffit que}
$\mu\geq\nu(D,\varphi)$.\qed

\medskip
Soit $O_L$ l'anneau des entiers de $L$, et soit $\pi\in O_F$ une
uniformisante qui est donc aussi une uniformisante de l'anneau de
valuation discr\`ete $O_L$. Soit $(D,\varphi)$ un isocristal de
dimension $d$. Un $O_L$-r\'eseau $M$ de $D$ est dit de type
$\mu\in
(\Z^d)_+$ s'il existe une base $e_1,\ldots, e_d$ de $M$ tel que
$\pi^{\mu_1}e_1,\ldots, \pi^{\mu_d}e_d$ soit une base de
$\varphi(M)$.
Tout $O_L$-r\'eseau $M$ a un type bien d\'etermin\'e que l'on note
$\mu(M)$.

\medskip\noindent
{\bf Th\'eor\`eme 2.} {\it Supposons $k$ alg\'{e}briquement clos.
Soit $(D,\varphi)$ un isocristal de
dimension $d$. Soit $\mu\in(\Z^d)_+$. Alors il existe un $O_L$-r\'eseau de
type $\mu$ dans $D$ si et seulement si $\mu\geq
\nu(D,\varphi)$.}\qed

\medskip
D'apr\`{e}s [F], Prop.\ 4.3.3, la condition d'admissibilit\'{e}
\'{e}quivaut
\`{a} demander  que
$P(\mu({\Cal F}^{\bullet}))$ et $P(\nu(D, \varphi))$
ont m\^{e}mes extr\'{e}mit\'{e}s et que, pour tout
sous-isocristal $(D',\varphi')$ de $(D,\varphi)$, si ${\Cal
F}^{\prime\bullet}$ est la filtration induite par ${\Cal
F}^{\bullet}$, alors
$P(\mu({\Cal F}^{\prime\bullet}))$ est en dessous de $P(\nu(D',
\varphi'))$.
L'implication directe du th\'eor\`eme 1 en r\'{e}sulte.

L'implication directe du th\'eor\`eme 2 est
l'in\'egalit\'e de Mazur ([Ka], Th.\ 1.4.1). Le contenu de ces deux
th\'eor\`emes
est
donc la r\'eciproque \`a ces in\'egalit\'es. Le lien entre eux
est donn\'e par un r\'{e}sultat de Laffaille [L], comme
on va le voir au \S 1.

Le th\'eor\`eme 2 est aussi obtenu, avec une preuve
diff\'erente, dans [KR]. En fait, on a une version de l'in\'egalit\'e de
Mazur pour un groupe r\'eductif quasi-d\'eploy\'e dans $F$ et
d\'eploy\'e
dans une extension non ramifi\'ee de $F$ [RR]. On peut conjecturer
que la
r\'eciproque \`a cette in\'egalit\'e est vraie, dans ce contexte
(existence de certains sous-groupes parahoriques
hypersp\'eciaux), et on
peut la d\'emontrer dans certains cas ([R], [KR]). En ce qui
concerne la
g\'en\'eralisation du th\'eor\`eme 1 dans cette direction, on a le
r\'esultat suivant~:

Soit $G$ un groupe r\'{e}ductif connexe {\it quasi-d\'eploy\'e} sur $F$.
Soit $A$ un tore d\'eploy\'e maximal. Soit $T$ le tore maximal contenant
$A$. Fixons un sous-groupe de Borel $B$ contenant $T$. Soit $\bar C$ la
chambre de Weyl ferm\'ee dans $X_*(T)\otimes \R$ correspondante. Si
$\bar F$ est une cl\^{o}ture alg\'{e}brique de $F$, le groupe
$\Gamma=\hbox{Gal}(\bar F/F)$ agit sur $\bar C$. Soit $\Afr=X_*(A)\otimes
\Q$ et $\Afr_+=\Afr\cap \bar C =\bar C^{\Gamma}\cap (X_*(T)\otimes \Q)$.

Soient $L'$ une extension finie de $L$ et $\mu:\G_m\to G$ un
morphisme d\'{e}fini sur $L'$. On note $\mu_0$ l'unique conjugu\'{e}
de $\mu$ qui appartient \`{a} $\bar C\cap X_*(T)$. On pose
$$
\bar\mu ={1\over (\Gamma:\Gamma_{\mu_0})}\cdot\sum\limits_{\tau\in
\Gamma/\Gamma_{\mu_0}}\tau\mu_0\in\Afr_+\ \ .\leqno(0.2)$$
C'est un \'{e}l\'{e}ment de $\Afr_+$ qui ne d\'{e}pend que de la
classe de conjugaison $\{\mu\}$ de $\mu$.

Soit $b\in G(L)$. On note $\bar\nu_b\in\Afr_+$ son point de Newton
([K3], introduction et \S 3.2). Rappelons comment on peut le
d\'{e}finir~: Soit $\Dbf$ le groupe diagonalisable sur $F$ de groupe
des caract\`{e}res $\Q$. On a $\hbox{Hom}(\Dbf,T)=\Afr$. A toute
repr\'{e}sentation rationnelle $(V,\varrho)$ de $G$ est associ\'{e}e
un isocristal $(D,\varphi)$, avec
$D=V\otimes_FL$
et $\varphi=\varrho(b)\cdot (id_V\otimes\sigma)$. La d\'{e}composition
par les pentes
$$D=\oplus_{\alpha\in\Q}D_{\alpha}$$
peut \^{e}tre consid\'{e}r\'{e}e comme un morphisme
$~\nu_{\varrho}:\Dbf\to GL(V)$ d\'{e}fini sur $L$.  Il existe un
unique morphisme $\nu_b:\Dbf\to G$ tel que
$\nu_{\varrho}=\varrho\circ\nu_b$ pour tout $\rho$. Alors $\bar\nu_b$
est l'unique conjugu\'{e} de $\nu_{b}$ sous l'action de $G(L)$
qui est dans $\Afr_+$. Il ne d\'{e}pend que de la classe de
$\sigma$-conjugaison de $b$ [K3].

   On a par ailleurs sur $\Afr_+$ l'ordre partiel usuel, pour lequel
   $\lambda\leq \lambda'$ si et seulement si $\lambda'-\lambda$ est
   combinaison lin\'eaire \`a coefficients $\geq 0$ des coracines
  simples relatives \`{a} $\bar C$.
\smallskip
Une paire $(b,\mu)$ form\'{e}e d'un \'el\'ement $b\in G(L)$ et d'un
homomorphisme $\mu:\G_m\to G$  d\'{e}fini sur une extension finie
$L'$ de $L$ est dite
{\it admissible} (appel\'ee {\it faiblement admissible} dans [RZ], d\'ef.
1.18)
    si, pour toute repr\'esentation rationnelle
$(V,\varrho)$ de $G$, l'isocristal associ\'{e} $(D,\varphi)$, muni de la
filtration
${\Cal F}^\bullet$ sur $D\otimes_LL'$ induite par
$\varrho\circ\mu$,
est admissible. Il suffit pour v\'erifier cette propri\'et\'e
de la tester sur une repr\'esentation fid\`ele.

\medskip\noindent
{\bf Th\'eor\`eme 3.} {\it Supposons $k$ alg\'{e}briquement clos.
Soit $G$  un groupe r\'{e}ductif connexe quasi-d\'eploy\'e sur $F$. Soient
$b\in
G(L)$, $L'$ une extension finie de $L$ et  $\mu : \G_m\to G$ un
morphisme d\'{e}fini sur $L'$. Pour qu'il
existe $\mu'\in\{ \mu\}$ d\'efini sur $L'$ et tel que la paire
$(b,\mu')$ soit admissible il faut et il suffit que
$\bar\mu\geq
\bar\nu_b$.}\qed

\medskip
Explicitons ce th\'eor\`eme dans le cas o\`u $G=T$ est un tore.
Dans ce
cas $\{\mu\}$ correspond \`a un seul \'el\'ement $\mu\in X_*(T)$ et
$\bar\mu$ est la moyenne sur l'orbite $\Gamma\cdot\mu$. Dans
ce cas
l'ordre partiel est trivial et l'\'{e}nonc\'{e} signifie que
$$(b,\mu)\ \hbox{est admissible}\ \Longleftrightarrow
\bar\mu
=\bar\nu_b\ \ .\leqno(0.3)$$ Ceci est exactement l'\'equivalence
entre (i)
et (iii) de la proposition 1.21 de [RZ]. Le th\'eor\`eme 3 est donc \`a
la fois
une g\'en\'eralisation de cette proposition, qui est le cas d'un tore,
et du th\'eor\`eme 1, qui est le cas de $GL_d$.

\medskip\noindent
{\bf Remarque :} Dans le th\'eor\`eme 3, soit $E\subset L'$ le corps de
d\'efinition de la classe de conjugaison $\{ \mu\}$ de $\mu$. Alors $E$
est une extension finie de $F$ et $L'$ contient $EL$. Invers\'ement,
\'etant donn\'ees une classe de conjugaison $\{\mu\}$ d\'efinie sur $E$ et
une extension $L'$ de $EL$, il existe $\mu\in \{\mu\}$ d\'efini sur $L'$,
cf.\ [K1], 1.4.3.

\medskip
Nous remercions G.\ Laumon pour des discussions utiles.

Durant la pr\'eparation de ce travail le deuxi\`eme auteur a
b\'en\'efici\'e du soutien financier du Minist\`ere de la Recherche et
de
la Fondation A. von Humboldt (prix Gay-Lussac / Humboldt), et
aussi de
l'hospitalit\'e des Universit\'es de Paris (Jussieu et Orsay).

\bigskip
\noindent {\bf 1 -- D\'emonstrations des th\'eor\`emes 1 et 2}
\medskip\noindent
Comme on l'a dit dans l'introduction, il s'agit seulement de
v\'{e}rifier
les implications r\'eciproques. On fixe l'isocristal $(D,\varphi)$
de
dimension $d$.

Soit $s\in\N$ tel que $s\cdot\nu (D,\varphi)\in\Z^d$. Soit $F'=F_s$
l'extension de degr\'e $s$ de $F$ contenue dans $L$. Pour tout
nombre
rationnel $\alpha$ tel que $s\alpha\in\Z$, soit $$ V_{\alpha}=\{ v\in
D;\
\varphi^s(v)= \pi^{s\alpha}v\}\ \ .\leqno(1.1)$$ Soit $V=\Sigma
V_{\alpha}$. Alors on sait que $V$ est somme directe des
$V_{\alpha}$ et
que l'application naturelle $V\otimes_{F'}L\to D$ est un
isomorphisme. En
plus, tout sous-isocristal $(D', \varphi')$ de $(D, \varphi)$ est
rationnel sur $F'$, i.e., il existe un unique sous-$F'$-espace
vectoriel
$V'$ de $V$ tel que $D'=V'\otimes_{F'}L$.

On fixe $\mu\in(\Z^d)_+$. Les filtrations ${\Cal F}^\bullet$ de type
$\mu$
sur $V\otimes_{F'}K$, pour $K$ extension de $F'$ variable, forment
l'ensemble des $K$-points d'une vari\'et\'e alg\'ebrique projective
$\hbox{\fett{$\Cal F$}}=\hbox{\fett{$\Cal F$}}(V,\mu)$ sur $F'$. En
fait,
$\hbox{\fett{$\Cal F$}}$ est une vari\'et\'e de drapeaux partiels de
$V$
et est de la forme $$\hbox{\fett{$\Cal F$}} =G/P\ \ ,\leqno(1.2)$$
o\`u
$G=GL(V)$ et o\`u $P$ est un sous-groupe parabolique de $G$.

Un point $x\in \hbox{\fett{$\Cal F$}}(L)$ sera dit {\it Weil-g\'en\'erique
relativement \`a} $F'$ si le point image du compos\'e
$${\roman{Spec}}\
L\buildrel x\over \longrightarrow \hbox{\fett{$\Cal
F$}}\times_{\roman{Spec}\ F'}{\roman{Spec}}\ L\longrightarrow
\hbox{\fett{$\Cal F$}}$$ est le point g\'en\'erique de
$\hbox{\fett{$\Cal
F$}}$.

\medskip\noindent
{\bf Lemme 1.1.} {\it (i) Il existe toujours des points $x\in
\hbox{\fett{$\Cal F$}} (L)$ Weil-g\'en\'eriques relativement \`{a} $F'$.
\par\noindent
(ii) Un tel point n'est contenu dans aucune sous-vari\'et\'e propre de
$\hbox{\fett{$\Cal F$}}$
    d\'efinie sur $F'$.}

\medskip\noindent
{\bf D\'emonstration.} (i) Par le lemme de Bruhat, on sait que le
corps des
fonctions $K$ de $\hbox{\fett{$\Cal F$}}$ est une extension
transcendante
{\it pure} de $F'$. Comme $L$ est de degr\'e de transcendance
infinie sur
$F'$, $K$ peut \^etre plong\'e dans $L$.
\par\noindent
(ii) C'est \'{e}vident.\qed

\medskip\noindent
{\bf Lemme 1.2.} {\it Soit $\hbox{\fett{$\Cal F$}}^{\bullet}$ une
filtration de $D$ de type $\mu$ correspondant \`a un point
Weil-g\'en\'erique relativement \`{a} $F'$ de $\hbox{\fett{$\Cal
F$}}(L)$. Soit $(D',
\varphi')$ un sous-isocristal de $(D,\varphi)$. Alors le type de la
filtration $\hbox{\fett{$\Cal F$}}^{\prime\bullet}$ induite par}
$\hbox{\fett{$\Cal F$}}^{\bullet}$ sur $D'$ est donn\'e par $$\mu({\Cal
F}^{\prime\bullet})= (\mu_{d-d'+1},\ldots,\mu_{d-1},\mu_d)\ \ ,\ \
\hbox{o\`u}\
d'={\roman{dim}}\ D'\ \ .$$
\medskip\noindent
{\bf D\'emonstration.} La filtration ${\Cal F}^{\bullet}$ est transverse
\`a $D'$, puisque les filtrations qui ne le sont pas forment une
sous-vari\'et\'e propre de $\hbox{\fett{$\Cal F$}}$ d\'efinie sur $F'$ (le
compl\'ement de l'unique orbite ouverte du parabolique standard de
type
$(d', d-d')$ sur $G/P$). On a donc $${\roman{dim}} ({\Cal F}^i\cap D')
=\max (0,{\roman{dim}}\ {\Cal F}^i -(d-d')),\ \forall i\in\Z\leqno(1.3)$$
et ceci implique l'assertion.\qed

\medskip\noindent
{\bf D\'emonstration du th\'eor\`eme 1.} Soient donc $\mu\geq
\nu(D,\varphi)$ et ${\Cal F}^\bullet$ une filtration de $D$
correspondant \`a un point Weil-g\'en\'erique relativement \`{a} $F'$.
Montrons
que ${\Cal F}^\bullet$ est admissible. La condition (i) de
l'admissibilit\'e $\vert\mu({\Cal F}^\bullet)\vert =
\vert\nu(D,\varphi)\vert$ (cf.\ introduction) r\'esulte de la d\'efinition
m\^eme de l'ordre partiel sur $(\Q^d)_+$. Soit $(D', \varphi')$ un
sous-isocristal de $(D,\varphi)$ de dimension $d'$. Alors
$(D',\varphi')$ est un facteur direct de $(D,\varphi)$ et son vecteur
des
pentes est de la forme $$\nu(D',\varphi')=(\nu_{j_1},\ldots,
\nu_{j_{d'}})\ \hbox{pour}\ 1\leq j_1<j_2<\ldots < j_{d'}\leq d\ \
.\leqno(1.4)$$ D'apr\`es le lemme 1.2, la filtration ${\Cal
F}^{\prime\bullet}$ sur $D'$ induite par ${\Cal F}^\bullet$ est de type
$\mu({\Cal F}^{\prime\bullet})= (\mu_{d-d'+1},\ldots , \mu_d)$.
L'in\'egalit\'e dans (ii) de la d\'efinition de l'admissibilit\'e
r\'esulte alors de $$\sum\limits_{i=1}^{d'}\mu_{d-d'+i}\leq
\sum\limits_{i=1}^{d'} \nu_{d-d'+i} \leq \sum\limits_{i=1}^{d'} \nu_{j_i}\
\ .\ \ \qed\leqno(1.5)$$

\medskip
Le th\'{e}or\`{e}me 2 r\'{e}sulte alors de la proposition suivante
(pour laquelle il n'est pas n\'{e}cessaire de supposer $k$
alg\'{e}briquement clos)~:
\medskip
\noindent
{\bf Proposition 1.3.} {\it Soient $(D,\varphi)$ un isocristal de
dimension $d$ et $\mu\in (\Z^{d})_+$.
Les assertions suivantes sont \'{e}quivalentes~:

i) Il existe une filtration admissible ${\cal F}^{\bullet}$ de type
$\mu$ sur $D$,

ii) il existe un $O_L$-r\'{e}seau de type $\mu$.}

\medskip\noindent
{\bf D\'emonstration.} C'est une cons\'{e}quence
du r\'esultat de Laffaille [L]~: Soit
$(D,\varphi)$ un isocristal et soit ${\Cal F}^\bullet$ une filtration
de $D$. Rappelons [L] qu'un $O_L$-r\'eseau $M$ dans $D$ est dit
{\it
adapt\'e \`a la filtration} ${\Cal F}^\bullet$ (ou, dans une autre
terminologie [FL], que le r\'eseau $M$, muni de sa filtration ${\Cal
F}^\bullet \cap M$, est {\it fortement divisible}), si l'on a
$$M=\varphi(\sum_i\pi^{-i} ({\Cal F}^i\cap M))\ \ .\leqno(1.6)$$
D'apr\`es le th\'{e}or\`{e}me 3.2 de
[L]\footnote{{{\it Stricto sensu}, ce r\'{e}sultat n'est \'{e}nonc\'{e}
dans [L] que lorsque   $F=\Q_p$ et la filtration est positive. Par
une torsion \`{a} la Tate, on se ram\`{e}ne au cas d'une filtration
positive. La preuve de Laffaille s'\'{e}tend alors
telle quelle au cas consid\'{e}r\'{e} ici~: il suffit de remplacer $p$ par
$\pi$.}\hfill}, la filtration ${\Cal F}^\bullet$ est faiblement
admissible
si et seulement s'il existe un r\'eseau $M$ adapt\'e \`a ${\Cal
F}^\bullet$. La proposition  r\'{e}sulte alors du lemme suivant~:

\medskip\noindent
{\bf Lemme 1.4.} {\it Soit $(D,\varphi)$ un isocristal de
dimension
$d$ et soit ${\Cal F}^\bullet$ une filtration de type $\mu({\Cal
F}^\bullet)$ sur $D$. Soit $M$ un r\'eseau adapt\'ee \`a ${\Cal
F}^\bullet$. Alors $$\mu(M)=\mu({\Cal F}^\bullet)\ \ .$$}
{\bf D\'{e}monstration.} Soit $u_1,u_2,\ldots,u_d$ une base de $M$
sur
$O_L$ {\it adapt\'{e}e \`{a} la filtration}, i.e. telle que, si, pour
$1\leq r\leq d$, $\mu_r$ d\'{e}signe le plus grand des entiers $i$
v\'{e}rifiant $u_r\in {\cal F}^{i}\cap M$, alors $${\cal F}^{i}\cap
M=\oplus_{\mu_r\geq i}O_Lu_r\ \ .$$ Quitte \`{a} changer l'ordre des
$u_r$, on peut supposer que
$\mu_1\geq\mu_2\geq\ldots\geq\mu_d$ et on
a alors $\mu({\cal F}^{\bullet})=(\mu_1,\mu_2,\ldots,\mu_d)$.

Si l'on pose $e_r=\pi^{-\mu_r}\varphi(u_r)$, on voit que
$e_1,e_2,\ldots,e_d$ est une base de $M$ sur $O_L$ tandis que
$\pi^{\mu_1}e_1,\pi^{\mu_2}e_2,\ldots,\pi^{\mu_d}e_d$ est une base
de $\varphi(M)$. On a donc $\mu(M)=(\mu_1,\mu_2,\ldots,\mu_d)$.
\qed
\medskip\noindent
{\bf Remarque :} Les \'{e}nonc\'{e}s des th\'{e}or\`{e}mes 1 et 2 ne
s'\'{e}tendent pas au cas o\`{u} le corps $k$ n'est plus
alg\'{e}briquement clos. Un exemple est donn\'{e} par $L=F$,
$(D,\varphi)=(L^{2},\roman{id})$, $\mu=(r,-r)$ o\`{u} $r$ est un entier
$\geq 1$. Alors $\mu\geq (0,0)=\nu(D,\varphi)$, mais il n'existe pas
de filtration admissible de type $\mu$.
\bigskip
\noindent {\bf 2 -- D\'emonstration du th\'eor\`eme 3}
\medskip\noindent
Soient $G$ quasi-d\'eploy\'e sur $F$, $L'$ une extension finie de $L$
et $\mu:\G_m\to G$ un morphisme d\'{e}fini sur $L'$. Avec les
notations de l'introduction, \`{a} $\mu$ correspondent des
\'{e}l\'{e}ments $\mu_0\in\bar C$ et $\bar\mu\in\Afr_+$ qui ne
d\'{e}pendent que de la classe de conjugaison $\{\mu\}$ de $\mu$.

Soit $(V,\varrho)$ une r\'epr\'esentation
rationnelle de $G$, o\`u $V$ est un espace vectoriel de dimension
$d$ sur
$F$. Alors on associe  \`a $\{\mu\}$ une
classe de conjugaison $\{ \varrho\circ\mu\}$ de morphismes de $\G_m$
dans  $GL(V)$ \`{a} laquelle correspond  un
\'el\'ement bien defini $\varrho(\mu_0)$ de $(\Z^d)_+$. De m\^eme,
$\bar\mu$ et par ailleurs n'importe quel \'el\'ement
$\bar\nu\in\Afr_+$
d\'efinissent des \'el\'ements $\varrho(\bar\mu)$ et $\varrho(\bar\nu)$ de
$(\Q^d)_+$.

\medskip\noindent
{\bf Lemme 2.1.}   {\it Soit $\bar\nu\in\Afr_+$.
Les
conditions suivantes sont \'equivalentes:

(i) $\bar\mu\leq\bar\nu$ (dans l'ordre partiel sur $\Afr_+$)

(ii) $\varrho(\bar\mu)\leq \varrho(\bar\nu)$ pour toute repr\'esentation
rationnelle $(V,\varrho)$ de $G$

(ii') $\varrho (\mu_0)\leq \varrho(\bar\nu)$ pour toute repr\'esentation
rationnelle $(V,\varrho)$ de $G$

(iii) comme en (ii), pour une repr\'esentation fid\`ele $(V,\varrho)$ de
$G$

(iii') comme en (ii'),  pour une repr\'esentation fid\`ele $(V,\varrho)$
de $G$}
\medskip\noindent
{\bf D\'emonstration.} On remarque que la classe de
conjugaison
$\{ \varrho\circ \mu\}$ est d\'efinie sur $F$, de sorte que
$\varrho(\mu_0)=\varrho(\bar\mu)$. Ceci \'etant, l'assertion r\'esulte
de
[RR], \S 2.\qed
\medskip\noindent
Ce lemme entra\^{\i}ne d\'eja l'implication directe dans le th\'eor\`eme
3, comme
cons\'e\-quence de l'implication directe dans le th\'eor\`eme 1.

Un cocaract\`{e}re $f : \G_m\to G$ d\'{e}finit, sur
chaque repr\'{e}sentation lin\'{e}aire $V$ de $G$, une graduation
$V=\oplus_{n\in\Z}V_n$  donc aussi une filtration
d\'{e}croissante en posant $\hbox{Fil}^{r}V=\oplus_{n\geq r}V_n$.
Deux cocaract\`eres de $G$ sont dit {\it par-\'equivalents} s'ils
d\'efinissent les m\^emes filtrations sur les repr\'{e}sentations
lin\'{e}aires de $V$. Les
classes
d'\'equivalence dans la classe de conjugaison $\{\mu\}$ de $\mu$
forment une vari\'et\'e
projective $\hbox{\fett{$\Cal F$}}=\hbox{\fett{$\Cal F$}}(G, \{\mu\})$
d\'efinie sur $L'$. On peut \'ecrire $$\hbox{\fett{$\Cal
F$}}(G,\{\mu\})=
G_{L'}/P_{\mu}\ \ ,\leqno(2.1)$$ o\`u $P_{\mu}$ est un sous-groupe
parabolique d\'efini sur $L'$.

Soit $b\in G(L)$. Pour
d\'emontrer le th\'eor\`eme 3 on peut remplacer $b\in G(L)$ par un
\'el\'ement $\sigma$-conjugu\'e $b'=gb\sigma(g)^{-1}$, pour $g\in
G(L)$,
${\roman{car}}(b,\mu)$ est admissible si et seulement si
$(gb\sigma(g)^{-1}, \ g\mu
g^{-1})$ l'est. D'apr\`es Kottwitz ([K2], \S 4) on peut donc supposer que
$b$
v\'erifie une identit\'e de d\'ecence ([RZ], p.~8), c'est-\`{a}-dire
qu'il existe un entier $s>0$ tel que $$(b\sigma)^s= (s\cdot
\nu_b)(\pi)\sigma^s$$
o\`{u} $\nu_b:\Dbf\to G$ a \'{e}t\'{e} d\'{e}fini avant l'\'{e}nonc\'{e}
du th\'{e}or\`{e}me 3 et $s$ d\'{e}signe l'endomor\-phisme de $\Dbf$
induit par la multiplication par $s$ sur $\Q$.

A $x\in{\Cal F}(L')$ on associe son {\it co-caract\`ere canonique de
Harder-Narasimhan} $$\lambda_x: \Dbf \longrightarrow
G\leqno(2.2)$$ qui
est bien d\'efini (sur $L'$) \`a par-\'equivalence pr\`es. Il est
caract\'eris\'e par
le fait que pour toute repr\'esentation rationnelle $(V,\varrho)$ de
$G$ la
$\Q$-filtration sur $V\otimes_FL$ induite par
$\varrho\circ\lambda_x$ est
la filtration canonique (index\'ee par les $HN$-pentes) de
l'isocristal filtr\'e $(V\otimes_FL,\ \varrho(b)\cdot ({\roman{id}}\otimes
\sigma), {\Cal F}^\bullet_{\varrho(x)})$, (cf. par exemple, [RZ], Prop.\
1.4). Soit $F'$ l'unique extension de degr\'{e} $s$ de $F$ contenue dans $L$.
On sait ([RZ], Prop.\ 1.36), que $\lambda_x$ est
d\'efini sur
$F'$.

Soit $G_{\roman{ab}} =G/G_{\roman{der}}$ le tore quotient maximal
de $G$.
L'isomorphisme $$\hbox{\fett{$\Cal F$}}
(G,\{\mu\})\buildrel\sim\over\longrightarrow \hbox{\fett{$\Cal F$}}
(G_{\roman{ad}}, \{\mu_{\roman{ad}}\})\times \hbox{\fett{$\Cal F$}}
(G_{\roman{ab}}, \{ \mu_{\roman{ab}}\})\leqno(2.3)$$ induit une
bijection
entre points (faiblement) admissibles, $$\hbox{\fett{$\Cal F$}}
(G,\{\mu\})(L')^{\roman{adm}}\buildrel\sim\over\longrightarrow
\hbox{\fett{$\Cal F$}}(G_{\roman{ad}}, \{\mu_{\roman{ad}}\})
(L')^{\roman{adm}}\times \hbox{\fett{$\Cal F$}} (G_{\roman{ab}},
\{\mu_{\roman{ab}}\})(L')^{\roman{adm}}\ \ .\leqno(2.4)$$ Pour
d\'emontrer le th\'eor\`eme 3, on peut donc supposer que $G$ est
semi-simple, puisque le cas d'un tore est d\'eja regl\'e par [RZ],
Prop.\
1.21, comme on l'a expliqu\'e dans l'introduction.

Soit donc $G$ semi-simple et soit $E'$ une extension finie de $F'$
contenue dans $L'$ sur laquelle la classe de conjugaison de $\mu$ est
d\'{e}finie. Soit $x\in \hbox{\fett{$\Cal F$}}
(L')$
Weil-g\'en\'erique relativement \`a $E'$. De tels points existent
d'apr\`es le
lemme
1.1 (lemme de Bruhat pour $G$). Montrons que, sous l'hypoth\`ese
$\bar\mu\geq \bar\nu_b$, $x$ est faiblement admissible. Il suffit de
voir
que le cocaract\`ere $\lambda_x$ associ\'e \`a $x$ est trivial.
Raisonnons
par l'absurde et supposons que le parabolique
$P_x=P_{\lambda_x}$ soit
propre. Ce parabolique est d\'efini sur $F'$.

\medskip\noindent
{\bf Lemme 2.2.} {\it Soient $G$ un groupe quasi-d\'eploy\'e sur $F$
et
$B$ un sous-groupe de Borel d\'efini sur $F$. Soit $P$ un sous-groupe
parabolique propre d\'efini sur $F'$.

(i) Il existe une repr\'esentation rationnelle irr\'eductible $(V,
\varrho)$ de $G$ et une droite $\ell\subset V\otimes_FF'$, telles
que
$P={\roman{Stab}}_G(\ell)$.

(ii) Pour une telle repr\'esentation soit $W_\bullet=((0)\subset
W_1\subset W_2\subset \ldots\subset W_d=V)$ un drapeau
complet de $V$
stable par $\varrho(B)$. Alors il existe $g\in G(F')$ tel que $\ell$ soit
transverse \`a $\varrho(g)\cdot W_\bullet$, i.e., tel que
$\ell\not\subset
\varrho(g)W_{d-1}$. }

\medskip\noindent
{\bf D\'emonstration.} (i) On peut \'evidemment supposer que $P$
est
standard, i.e., contient $B\otimes_F F'$. Soit $\bar C^*\subset
X^*(T)\otimes\R$
la
chambre de Weyl ferm\'ee dans l'espace des caract\`eres
correspondant \`a
$B$. Alors $P$ correspond \`a une facette de $\bar C^*$. Soit
$\lambda\in
X^*(T)$ \`a l'int\'erieur de cette facette, et soit $V_{\lambda}$ la
repr\'esentation irr\'eductible de plus haut poids $\lambda$. Alors
$V_{\lambda}$ est un espace vectoriel sur $F'$ sur lequel agit
$G_{F'}$. Soit $F_\lambda=F^{\Gamma_\lambda}$. Alors
$V_\lambda$ est
d\'efini sur $F_\lambda$, i.e. est de la forme
$V_\lambda=V_\lambda^0\otimes_{F_\lambda}{F'}$ et
$V_{\lambda}^0$ est
irr\'eductible comme repr\'esentation de $G_{F_\lambda}$. Le $F$-
espace
vectoriel cherch\'e $V$ est \'egale \`a $V_{\lambda}^0$ et est muni
de la
repr\'esentation rationnelle $\varrho$ de $G$ d\'efinie par la
composition des
homomorphismes canoniques,
$$G\to R_{F_\lambda/F}(G_{F_\lambda})\to
R_{F_\lambda/F}(GL_{F_{\lambda}}(V_{\lambda}^0))=
GL_F(V),\leqno(2.5)$$
(cf.\ [T], th. 7.2).
Alors $(V,\varrho)$
r\'epond
\`a la question en prenant pour $\ell$ la droite engendr\'e par un
vector
de poids $\lambda$ dans $V_{\lambda}$.

(ii) Le sous-espace de $V\otimes_FF'$ engendr\'e par $\{
\varrho(g)x;\
g\in G(F'),\ x\in\ell\}$ est une repr\'esentation irr\'eductible de
$G_{F'}$. Le plus petit sous-espace d\'efini sur $F$ et le contenant
est
une repr\'esentation de $G$. Son espace ne peut donc pas \^etre
contenu
dans $W_{d-1}$.\qed

\medskip
Appliquons ce lemme au parabolique $P_{\lambda_x}$. On
consid\`ere l'isocristal $(D,\varphi)=(V\otimes_FL,\varrho(b)\cdot
({\roman{id}}\otimes\sigma))$, muni de sa filtration ${\Cal
F}^\bullet_{\varrho(x)}$ sur $D\otimes_LL'$. Soit
$D'=\ell\otimes_{F'}L$. Alors $D'$ est un cran dans la filtration de
Harder-Narasimhan de $D$, car il est stable par le parabolique de
Harder-Narasimhan $P_{\varrho\circ \lambda_x}$ de $GL(V)_L$.
Le type
$\mu({\Cal F}^{\prime\bullet})\in\Z$ de la filtration
induite par ${\Cal F}^\bullet_{\varrho(x)}$ sur
$D'\otimes_LL'$ et la pente $\nu(D',\varphi')\in\Q$ de
l'isocristal $(D', \varphi')= (D', \varphi\vert D')$ v\'{e}rifient donc
$$\mu({\Cal
F}^{\prime\bullet}) > \nu(D', \varphi')\ \ .\leqno(2.6)$$ Le lemme 2.2
implique que les points $y\in \hbox{\fett{$\Cal F$}}(L')$ tels
que
la filtration ${\Cal F}^\bullet_{\varrho(y)}$ de $D\otimes_LL'$
ne
soit pas transverse \`a $D'\otimes_LL'$ forment une sous-vari\'et\'e
ferm\'e {\it propre} de $\hbox{\fett{$\Cal F$}}$ d\'efini sur $E'$.
Comme
$x$ est Weil-g\'en\'erique relativement \`a $E'$, il ne peut pas \^etre
contenu
dans
cette sous-vari\'et\'e. Soit ${\Cal F}^\bullet_{\varrho(x)}$ de type
$(\mu_1,\ldots, \mu_d)\in (\Z^d)_+$. Comme dans le lemme 1.2, la
transversalit\'{e} implique que
   $$\mu({\Cal F}^{\prime\bullet})=\mu_d\ \
.\leqno(2.7)$$
Soit $\nu(D,\varphi)=(\nu_1,\ldots,\nu_d)\in (\Q^d)_+$. Alors
$\nu(D',\varphi')=\nu_j$, pour un $j$ avec $1\leq j\leq d$. Comme
$\bar\mu\geq \bar\nu_{\bar b}$, le lemme 2.1 nous montre que
$(\mu_1,\ldots,
\mu_d)\geq (\nu_1,\ldots,\nu_d)$, ce qui nous  donne $$\mu({\Cal
F}^{\prime\bullet})= \mu_d\leq \nu_d\leq \nu_j=\nu(D', \varphi')\ \
,\leqno(2.8)$$ contredisant (2.6).\qed

\bigskip
\noindent {\bf 3 -- Remarques suppl\'ementaires}
\medskip\noindent
Supposons $k$ alg\'{e}briquement clos.
\medskip\noindent

a) D'apr\`{e}s la proposition 1.3,  dans le cas de $GL_d$,
les deux th\'eor\`emes d'existence sont \'{e}quivalents. Il faut
pourtant
souligner que dans le contexte d'un groupe r\'eductif $G$,
quasi-d\'eploy\'e sur $F$ et d\'eploy\'e sur une extension non
ramifi\'ee de
de $F$, l'existence d'un sommet hyper-sp\'ecial d'un type $\mu$
implique,
outre l'in\'egalit\'e $\bar\mu\geq \bar\nu_b$ qui appara\^{\i}t  dans
le
th\'eor\`eme 3, l'identit\'e
$$\mu^{\natural}=\kappa(b)\leqno(3.1)$$ dans
$\pi_1(G)_{\Gamma}$ (cf.\ [RR], \S 4) . On voit
donc qu'a priori, dans le cas g\'{e}n\'{e}ral,
la question d'existence de r\'eseaux est plus
subtile que
la question d'existence de filtrations. Le cas de $GL_d$ est
exceptionnel car alors (3.1) r\'{e}sulte de l'in\'{e}galit\'{e}
    $\bar\mu\geq \bar\nu_b$. Le lien entre le th\'eor\`eme 1
et
le th\'eor\`eme 2 se fait via le th\'eor\`eme de Laffaille sur l'existence
de r\'eseaux fortement divisibles et on peut se demander s'il existe
une notion analogue en th\'eorie de groupes.

\medskip\noindent
b) La construction de Laffaille d'un r\'eseau adapt\'e \`a une filtration
admissible fournit \`a partir d'un r\'eseau $M$ arbitraire un unique
r\'eseau adapt\'e maximal $M^{\roman{max}}$ contenu dans $M$,
et aussi un
unique r\'eseau adapt\'e minimal $M^{\roman{min}}$ contenant $M$.
Malheureusement, cette construction n'est pas compatible
aux
operations d'alg\`ebre lin\'eaire usuelles (produit tensoriel:
$(M_1\otimes M_2)^{\roman{max}}\neq M_1^{\roman{max}}\otimes
M_2^{\roman{max}}$, $(M_1\otimes M_2)^{\roman{min}}\neq
M_1^{\roman{min}}\otimes M_2^{\roman{min}}$, passage au dual:
on a $(M^*)^{\roman{max}}= (M^{\roman{min}})^*$ et non
$(M^{\roman{max}})^*$). Ceci, tout comme la remarque
pr\'{e}c\'{e}dente, semble montrer qu'il n'est
pas possible de d\'eduire du th\'eor\`eme 2 et du
th\'eor\`eme 3
pour le groupe symplectique l'existence de r\'eseaux autoduaux
d'un type
donn\'e dans un isocristal symplectique (un tel th\'eor\`eme
d'existence est pourtant d\'emontr\'e, par d'autres m\'{e}thodes,
dans [KR]).

\medskip\noindent
c) Nous ignorons ce qui se passe si l'on abandonne l'hypoth\`ese
$\bar\mu\geq \bar\nu_b$. Explicitons dans le cas de $GL_d$. Soit
$(D,\varphi)$ un isocristal de dimension $d$. Soit $\mu\in
(\Z^d)_+$
avec $\vert\mu\vert=\vert\nu(D,\varphi)\vert$, et soit
$\hbox{\fett{$\Cal
F$}}$ la vari\'et\'e des filtrations de type $\mu$ sur $D$. Soit
$x\in\hbox{\fett{$\Cal F$}}(L)$ et soit $\lambda_x:\Dbf\to GL(D)$
le
$HN$-cocaract\`ere associ\'e (bien d\'efini \`a par-\'equivalence
pr\`es).
Alors $\lambda_x$ d\'efinit un \'{e}l\'{e}ment bien d\'etermin\'e
$\lambda(x)\in
(\Q^d)_+$ avec $\vert\lambda(x)\vert =0$. On obtient ainsi une
d\'ecomposition disjointe $$\hbox{\fett{$\Cal
F$}}(L)=\bigcup\limits^\bullet_{\lambda\in (\Q^d)_+,\atop
\vert\lambda\vert =0}\hbox{\fett{$\Cal F$}}(L)_{\lambda}\ \
,\leqno(3.2)$$
o\`u $$\hbox{\fett{$\Cal F$}}(L)_{\lambda}=\{ x\in \hbox{\fett{$\Cal
F$}}(L);\ \lambda(x)=\lambda\}\ \ .$$

Si $\mu\geq\nu(D,\varphi)$, le th\'eor\`eme 1 montre que l'unique
\'el\'ement minimal $\lambda=0$ de l'ensemble d'indices donne une
contribution non-vide \`a (3.2). M\^eme dans le cas $\mu\geq
\nu(D,\varphi)$ nous n'avons pas d\'etermin\'e l'ensemble des
indices dans
(3.2) avec une contribution non-vide. Si $\mu\not\geq
\nu(D,\varphi)$ nous
ignorons m\^eme les indices {\it minimaux} parmi ceux avec une
contribution non-vide \`a (3.2).

\medskip\noindent
d) Rappelons que $C_p$ d\'{e}signe le compl\'et\'e d'une
cl\^{o}ture alg\'ebrique de $L$. Soient $G,b,L',\mu$ comme dans
l'\'enonc\'e du
th\'eor\`eme 3. Alors on sait [RZ], Prop. 1.36 que les points
admissibles par rapport \`a $b$ dans $\hbox{\fett{$\Cal
F$}}(G,\{\mu\})(C_p)$ forment un espace rigide-analytique sur $L'$, ouvert
admissible de $\hbox{\fett{$\Cal F$}}(G,\{\mu\})$,
que
nous noterons $\hbox{\fett{$\Cal F$}}(G,\{\mu\})^{\roman{adm}}$.
On peut
alors reformuler ainsi le th\'eor\`eme 3~:

\medskip\noindent
{\bf Th\'eor\`eme 3'.} {\it On suppose $k$ alg\'{e}briquement clos.
Soit $G,b,L',\mu$ comme dans le
th\'eor\`eme 3. Les conditions sui\-vantes sont \'equivalentes:
\par\noindent
(i) $\hbox{\fett{$\Cal F$}}(G,\{\mu\})^{\roman{adm}}\neq\emptyset$
\par\noindent
(ii) $\hbox{\fett{$\Cal F$}}(G,\{\mu\})^{\roman{adm}}(L')\neq\emptyset$
\par\noindent
(iii) $\bar\mu\geq \bar\nu_b$. }

\medskip\noindent
En effet, (i) $\Rightarrow$ (iii) se d\'emontre comme dans le \S 2
par
r\'eduction \`a l'in\'egalit\'e [F] 4.3.3 pour le cas de $GL_d$.
L'implication (iii) $\Rightarrow$ (ii) est le th\'eor\`eme 3 et
l'implication (ii) $\Rightarrow$ (i) est triviale.\qed

\medskip\noindent
e) Soit $(D,\varphi,N)$ un $(\varphi,N)$-module sur $L$ (cf. par
exemple, [CF]),
c'est-\`{a}-dire un isocristal muni d'une application $L$-lin\'{e}aire
$N: D\rightarrow D$ telle que $N\varphi=q\varphi N$ (o\`{u} $q$,
rappelons-le, est le nombre d'\'{e}l\'{e}ments du corps r\'{e}siduel
de $F$). On a encore
une
notion de filtration admissible dans ce cadre et le th\'{e}or\`{e}me
1 s'\'{e}tend~: cela r\'{e}sulte de ce que toute filtration
admissible sur l'isocristal sous-jacent (obtenu en oubliant l'action
de $N$) est a fortiori admissible sur le $(\varphi,N)$-module.
\bigskip
\centerline{\bf Bibliographie}
\bigskip\noindent

\ref{[CF]} P.\ Colmez, J.-M.\ Fontaine: Construction des
repr\'{e}sen\-tations semi-stables. Invent.\ Math. {\bf 140} (2000), 1--
43
\ref{[De]} J.-M. Decauwert: Classification des $A$-modules formels,
C.R.A.S. Paris {\bf 282} (1976), 1413-1416
    \ref{[F]} J.-M.\ Fontaine: Modules galoisiens, modules filtr\'es et
anneaux de Barsotti-Tate. Ast\'erisque {\bf 65} (1979), 3--80
    \ref{[FL]} J.-M.\ Fontaine, G.\ Laffaille: Construction de
repr\'esentations
$p$-adiques.  Ann. Sci. \'Ecole Norm. Sup. {\bf 15} (1982), 547--608
    \ref{[Ka]} N.\ Katz: Slope filtration of $F$-crystals.
    Ast\'erisque {\bf 63} (1979), 113-163
    \ref{[K1]} R.\ Kottwitz: Shimura varieties and twisted orbital
integrals.
Math. Ann. {\bf 269} (1984), 287--300.
    \ref{[K2]} R.\ Kottwitz:
Isocrystals with additional structure. Comp.\ Math. {\bf 56} (1985),
201-222
    \ref{[K3]} R.\ Kottwitz: Isocrystals with additional structure.
II. Comp.\ Math.\ {\bf 109} (1997),  255-339
    \ref{[KR]} R.\ Kottwitz, M.\ Rapoport: On
the existence of $F$-crystals. In preparation
    \ref{[L]} G.\ Laffaille: Groupes $p$-divisibles et modules filtr\'es: le
cas
    peu ramifi\'e.  Bull.\ Soc.\ Math.\ France {\bf 108} (1980),  187--206
    \ref{[R]} M.\ Rapoport: A positivity property of the Satake
isomorphism. Manu\-scrip\-ta Math. {\bf 101} (2000), 153--166
     \ref{[RR]} M.\ Rapoport, M.\ Richartz:
     On the classification and specialization of $F$-isocrystals with
additional
structure. Comp.\ Math. {\bf 103} (1996), 153--181
    \ref{[RZ]} M.\
Rapoport, Th.\ Zink: Period spaces for $p$-divisible groups. Annals
of
Mathematics Studies {\bf 141}, Princeton University Press,
Princeton, NJ,
1996
    \ref{[T]} J.\ Tits: Repr\'esentations lin\'eaires irr\'eductibles d'un
    groupe r\'eductif sur un corps quelconque. J.\ Reine Angew. Math.
{\bf
    247} (1971), 196--220
\end